\documentclass[a4paper,10pt]{article}
\usepackage{amsfonts}
\usepackage{amstext}
\usepackage{amsthm}
\usepackage{amsmath}
\usepackage{amssymb}
\usepackage{latexsym}
\usepackage[latin1]{inputenc}
\newcommand{\bl}{\hfill\rule{2mm}{2mm}}
\newcommand{\R}{\Bbb{R}}

\newtheorem{teor}{Theorem}[section]
\newtheorem{propo}{Proposition}[section]

\newtheorem{lema}{Lemma}[section]

\newcommand{\n}{\noindent}

\begin{document}

\title{Sharp $L^p$-entropy inequalities on manifolds
 }

\author{\textbf{Jurandir Ceccon \footnote{\textit{E-mail addresses}:
ceccon@ufpr.br (J. Ceccon)}}\\ {\small\it Departamento
de Matem\'{a}tica, Universidade Federal do Paran\'{a},}\\ {\small\it Caixa Postal 19081, 81531-980, Curitiba, PR, Brazil}\\
\textbf{Marcos Montenegro \footnote{\textit{E-mail addresses}:
montene@mat.ufmg.br (M. Montenegro)}}\\ {\small\it Departamento de Matem\'{a}tica,
Universidade Federal de Minas Gerais,}\\ {\small\it Caixa Postal 702, 30123-970, Belo Horizonte, MG, Brazil}} \maketitle

\markboth{abstract}{abstract}
\addcontentsline{toc}{chapter}{abstract}


\begin{center}
\small{
{\bf Abstract}}
\end{center}

\small{In 2003, Del Pino and Dolbeault \cite{DPDo} and Gentil \cite{G} investigated, independently, best constants and extremals associated to sharp Euclidean $L^p$-entropy inequalities. In this work, we present some important advances in the Riemannian context. Namely, let $(M,g)$ be a compact Riemannian manifold of dimension $n \geq 3$. For $1 < p \leq 2$, we prove that the
sharp Riemannian $L^p$-entropy inequality

\[
\int_M |u|^p \log(|u|^p) dv_g \leq \frac{n}{p} \log \left( {\cal A}_{opt} \int_M |\nabla_g u|^p dv_g + {\cal B} \right)
\]

\n holds on all functions $u \in H^{1,p}(M)$ such that $||u||_{L^p(M)} = 1$. Moreover, we show that the first best Riemannian constant ${\cal A}_{opt}$ is equal to the corresponding Euclidean one. Our
approach is inspired on the Bakry, Coulhon, Ledoux and Sallof-Coste's idea \cite{Ba} of getting Euclidean entropy inequalities as a limit case of suitable Gagliardo-Nirenberg inequalities. It is conjectured that the above inequality sometimes fails for $p > 2$.}

\vspace{0,5cm}

\begin{center}
\section{Introduction}
\end{center}

\n 1.1. {\bf Overview} Logarithmic Sobolev inequalities are a powerful tool in Real Analysis, Complex Analysis, Geometric Analysis, Convex Geometry and Probability. The pioneer work by L. Gross \cite{Gr} put forward the equivalence between a class of logarithmic Sobolev inequalities and hypercontractivity of the associated heat semigroup. Particularly, the Gaussian logarithmic Sobolev inequality has important applications such as the behavior of the Perelman entropy functional along the Ricci flow, which it plays a key role in the proof program of the Poincar\'{e} Conjecture \cite{Pe}, the Talagrand transport-entropy inequality within the optimal transport theory \cite{Tal}, the concentration of measure in probability theory \cite{Le1}, among many others.

Recent developments has pointed out a close relationship between a class of sharp logarithmic Sobolev inequalities, known also as sharp entropy inequalities, and hypercontrativity for some nonlinear diffusion equations. We refer the reader to \cite{DDG}, \cite{G0} and \cite{G} for references in the Euclidean context and \cite{BGL}, \cite{BGL1} and \cite{G} in the Riemannian one.

The main interest here is best constants and sharp entropy inequalities within the Riemannian setting. Before going further and describing some problems of interest, we present the Euclidean corresponding ones as well as their solutions.

The Euclidean $L^p$-entropy inequality states that, for any function $u \in W^{1,p}(\R^n)$ with $\int_{\R^n} |u|^p dx = 1$,

\begin{equation}\label{dee}
\int_{\R^n} |u|^p \log(|u|^p) dx \leq \frac{n}{p} \log \left({\cal A}_0(p)
\int_{\R^n} |\nabla u|^p dx\right)\; ,
\end{equation}

\n where $n \geq 2$, $p \geq 1$ and ${\cal A}_0(p)$ denotes the best entropy constant.

In 1978, Weissler \cite{W} presented the sharp inequality (\ref{dee}) for $p = 2$, which it is equivalent to the Gaussian logarithmic Sobolev inequality due to Gross \cite{Gr}. Later, Carlen \cite{Carlen} showed that the extremals of (\ref{dee}) are precisely dilations and translations of the Gaussian

\[
u(x) = \pi^{-\frac n2} e^{-|x|^2}\; .
\]

\n When $p = 1$, Ledoux \cite{Le} established the inequality (\ref{dee}) and Beckner \cite{Be} classified its extremals as normalized characteristic functions of balls. In this same work, Beckner showed the validity of (\ref{dee}) for $1 < p < n$ and Del Pino and Dolbeault \cite{DPDo} proved that the extremals are precisely dilations and translations of the function

\[
u(x) = \pi^{-\frac n2} \frac{\Gamma(\frac{n}{2} + 1)}{\Gamma(\frac{n(p - 1)}{p} + 1)} e^{-|x|^{\frac{p}{p-1}}}\; .
\]

\n Finally, thanks to the uniqueness argument exhibited in \cite{DPDo} for some elliptic PDEs, Gentil \cite{G} established the validity of (\ref{dee}) and extended the above classification for any $p \geq n$. From this, they concluded for $p > 1$ that

\[
{\cal A}_0(p) = \frac{p}{n} \left( \frac{p - 1}{e} \right)^{p - 1} \pi^{-\frac{p}{2}} \left( \frac{\Gamma(\frac{n}{2} +
1)}{\Gamma(\frac{n(p - 1)}{p} + 1)} \right)^{\frac{p}{n}}\; .
\]

\n In particular,

\[
{\cal A}_0(2) = \frac{2}{n \pi e}\; .
\]

It deserves mention a key idea introduced by Bakry, Coulhon, Ledoux and Salof-Coste in \cite{Ba}. Namely, they discovered for $1 \leq p < n$ that the non-sharp Euclidean $L^p$-entropy inequality can be deduced as a limit case of certain non-sharp Euclidean Gagliardo-Nirenberg inequalities. Inspired on this remark, Del Pino and Dolbeault \cite{DPDo} considered a class of sharp Euclidean Gagliardo-Nirenberg inequalities which interpolates the sharp $L^p$-Sobolev and $L^p$-entropy inequalities for $1 < p < n$ and utilized symmetry results of positive solutions of elliptic PDEs in order to characterize all extremals of (\ref{dee}). In turn, the tools used by Gentil in \cite{G} consist of an equivalence between inequality (\ref{dee}) and hypercontractivity of the Hamilton-Jacobi semigroup and of the proof of hypercontractivity via the Pr\'{e}kopa-Leindler inequality, also known as generalized Brunn-Minkowski inequality.\\

\n 1.2. {\bf Riemannian setting} Before stating extensions, problems and results, a little bit of definition and notation should be introduced.

Let $(M,g)$ be a smooth compact Riemannian manifold of dimension $n \geq 2$ and $p \geq 1$ be a parameter. The {\bf head $L^p$-entropy inequality} states that there exist constants ${\cal A}, {\cal B} \in \R$ such that, for any $u \in H^{1,p}(M)$ with $\int_M |u|^p dv_g = 1$,

\begin{gather}\label{AB}
\int_{M} |u|^p \log |u|^p\; dv_g \leq \frac{n}{p} \log \left( {\cal A} \int_M |\nabla_g u|^p dv_g + {\cal B} \right)\; , \tag{$L({\cal A},{\cal B})$}
\end{gather}
where $dv_g$ denotes the Riemannian volume element, $\nabla_g$ is the gradient operator of $g$ and $H^{1,p}(M)$ is the Sobolev space defined as the completion of $C^{\infty}(M)$ under the norm

\[
||u||_{H^{1,p}(M)} := \left( \int_{M} |\nabla_g u|^p\; dv_g  + \int_{M} |u|^p\; dv_g \right)^{\frac{1}{p}}\; .
\]

The following definitions and notations related to (\ref{AB}) are quite natural when one desires to extend the sharp $L^p$-entropy inequality (\ref{dee}) to the compact Riemannian context.

The {\bf first Riemannian $L^p$-entropy best constant} is defined by

\[
{\cal A}_0(p,g) := \inf \{ {\cal A} \in \R: \mbox{ there exists} \hspace{0,18cm} {\cal B} \in \R \hspace{0,18cm} \mbox{such that (\ref{AB})} \hspace{0,18cm} \mbox{holds for all} \hspace{0,18cm} u \in H^{1,p}(M) \hspace{0,18cm} \mbox{with} \hspace{0,18cm} \int_M |u|^p dv_g = 1\}\; .
\]

The {\bf first sharp Riemannain $L^p$-entropy inequality} states that there exists a constant ${\cal B} \in \R$ such that, for any $u \in H^{1,p}(M)$ with $\int_M |u|^p dv_g = 1$,

\[
\int_{M} |u|^p \log |u|^p\; dv_g \leq \frac{n}{p} \log \left( {\cal A}_0(p,g) \int_M |\nabla_g u|^p dv_g + {\cal B} \right)\; .
\]

If the preceding inequality is true, then we can define the {\bf second Riemannian $L^p$-entropy best constant} as

\[
{\cal B}_0(p,g) := \inf \{ {\cal B} \in \R:\; (L({\cal A}_0(p,g),{\cal B})) \hspace{0,18cm} \mbox{holds for all} \hspace{0,18cm} u \in H^{1,p}(M) \hspace{0,18cm} \mbox{with} \hspace{0,18cm} \int_M |u|^p dv_g = 1\}
\]
and the {\bf second sharp Riemannian $L^p$-entropy inequality} as the saturated version of ($L({\cal A}_0(p,g),{\cal B})$) on the functions $u \in H^{1,p}(M)$ with $\int_M |u|^p dv_g = 1$, that is

\[
\int_{M} |u|^p \log |u|^p\; dv_g \leq \frac{n}{p} \log \left( {\cal A}_0(p,g) \int_M |\nabla_g u|^p dv_g + {\cal B}_0(p,g) \right)\; .
\]

\n Note that ($L({\cal A}_0(p,g),{\cal B}_0(p,g))$) is sharp with respect to both the first and second best constants in the sense that none of them can be lowered. In a natural way, it arises then the notion of extremals of ($L({\cal A}_0(p,g),{\cal B}_0(p,g))$). A function $u_0  \in H^{1,p}(M)$ satisfying $\int_M |u_0|^p dv_g = 1$ is said to be extremal, if

\[
\int_{M} |u_0|^p \log |u_0|^p\; dv_g = \frac{n}{p} \log \left( {\cal A}_0(p,g) \int_M |\nabla_g u_0|^p dv_g + {\cal B}_0(p,g) \right)\; .
\]

\n We denote the set of extremals of ($L({\cal A}_0(p,g),{\cal B}_0(p,g))$) by ${\cal E}_0(p,g)$.

As pointed out by Brouttelande in \cite{Bro}, the knowledge of values of ${\cal A}_0(2,g)$ and ${\cal B}_0(2,g)$ and the validity of ($L({\cal A}_0(2,g),{\cal B}_0(2,g))$) play a fundamental role in the study of hypercontractivity and/or ultracontractivity of the heat semigroup on compact manifolds. Indeed, let $(P_t)_{t \geq 0}$ be the heat semigroup related to a smooth compact Riemannian manifold $(M,g)$ of dimension $n \geq 2$. We recall the following result due to Bakry \cite{Ba1}:

\begin{teor} \label{Ba}
Assume that, for any $u \in H^{1,2}(M)$ with $\int_M u^2 dv_g = 1$,

\[
\int_{M} u^2 \log u^2\; dv_g \leq \phi \left( \int_M |\nabla_g u|^2 dv_g\right)\; ,
\]

\n where $\phi: \R_+ \rightarrow \R_+$ is a concave, increasing and of class $C^1$ function. Then, for any $1 \leq p < q \leq \infty$,

\[
||P_t f||_{L^q(M)} \leq e^m ||f||_{L^p(M)}\; ,
\]

\n where

\[
t = \int^q_p \frac{\phi'(v(s))}{4(s-1)}\; ds\ \ \mbox{and}\ \ m = \int^q_p \frac{\phi(v(s)) - v(s)\phi'(v(s))}{s^2}\; ds\; ,
\]

\n provided we find a function $v \geq 0$ for which these two integrals are finite.
\end{teor}

Assuming ${\cal A}_0(2,g) = {\cal A}_0(2)$, ($L({\cal A}_0(2),{\cal B}_0(2,g))$) is valid and choosing

\[
\phi(x) = \frac{n}{2} \log({\cal A}_0(2) x + {\cal B}_0(2,g))\ \ {\rm and}\ \ v(s) = \frac{\lambda s^2}{s-1} - \frac{{\cal B}_0(2,g)}{{\cal A}_0(2)}\; ,
\]

\n one deduces (see \cite{Bro}) that

\[
||P_t f||_{L^\infty(M)} \leq \frac{1}{(4 \pi t)^{n/2}} e^{\frac{2{\cal B}_0(2,g)}{3{\cal A}_0(2)}t} ||f||_{L^1(M)}
\]

\n for all $0 < t \leq \frac{n}{2} \frac{{\cal A}_0(2)}{{\cal B}_0(2,g)}$. In particular, this conclusion shows that upper and lower bounds of ${\cal B}_0(2,g)$ are important. Note that $||P_t||_{1,\infty}$ is of the order of $\frac{1}{(4 \pi t)^{\frac{n}{2}}}$ for $t > 0$ small enough. For comparison, we recall that the heat semigroup on $\R^n$ satisfies

\[
||P_t||_{1,\infty} = \frac{1}{(4 \pi t)^{\frac{n}{2}}}
\]

\n for any $t > 0$. Applications of this nature make us wonder about some issues related to the best constants ${\cal A}_0(p,g)$ and ${\cal B}_0(p,g)$ and to sharp inequality ($L({\cal A}_0(p,g),{\cal B}_0(p,g))$). A number of key questions make up the following Riemannian entropy program:

\begin{itemize}

\item[(a)] Is ${\cal A}_0(p,g)$ well defined?

\item[(b)] Is it possible to find exact values and/or bounds for ${\cal A}_0(p,g)$ and ${\cal B}_0(p,g)$?

\item[(c)] Is there a constant ${\cal B} \in \R$ such that ($L({\cal A}_0(p,g),{\cal B})$) is valid for all $u \in H^{1,p}(M)$ with $\int_M |u|^p dv_g = 1$?

\item[(d)] Do ${\cal A}_0(p,g)$ and ${\cal B}_0(p,g)$ depend continuously on $p$ and $g$ in some topology?

\item[(e)] Is ${\cal E}_0(p,g)$ non-empty and compact in some topology?

\end{itemize}

Perhaps, contrary to what one might expect, the issue (a) is not straightforward. The main difficult in Riemannian $L^p$-entropy inequalities is that local-to-global kind arguments do not work well. In other words, the normalization condition $\int_M |u|^p dv_g = 1$ and the involved log functions in (\ref{AB}) do not allow a comparison to the flat corresponding situation. In conclusion, it is not clear that (\ref{AB}) is true for some constants ${\cal A}$ and ${\cal B}$ and also that, for any $\varepsilon > 0$, there exists a constant ${\cal B}_\varepsilon \in \R$ such that

\[
\int_{M} |u|^p \log |u|^p\; dv_g \leq \frac{n}{p} \log \left( ({\cal A}_0(p) + \varepsilon) \int_M |\nabla_g u|^p dv_g + {\cal B}_\varepsilon \right)
\]

\n holds for all $u \in H^{1,p}(M)$ with $\int_M |u|^p dv_g = 1$.

Recent contributions on sharp Euclidean entropy inequalities, developments over thirty years in the field of sharp Riemannian Sobolev inequalities and results and remarks contained in this work suggest the following conjectures:

\begin{itemize}

\item[(I)] (well-posedness) ${\cal A}_0(p,g)$ is well defined for all $p \geq 1$

\item[(II)] (first best constant) ${\cal A}_0(p,g) = {\cal A}_0(p)$ for all $p \geq 1$

\item[(III)] (validity) For any $1 \leq p \leq 2$, there exists a constant ${\cal B} \in \R$ such that ($L({\cal A}_0(p,g),{\cal B})$) holds for all $u \in H^{1,p}(M)$ with $\int_M |u|^p dv_g = 1$

\item[(IV)] (non-validity) If $(M,g)$ has positive scalar curvature somewhere, then for any $p > 2$ and any constant ${\cal B} \in \R$, there exists $u_{{\cal B}} \in H^{1,p}(M)$ with $\int_M |u_{{\cal B}}|^p dv_g = 1$ such that ($L({\cal A}_0(p,g),{\cal B})$) fails for $u_{{\cal B}}$

\item[(V)] (extremal existence) ${\cal E}_0(p,g)$ is non-empty if either $1 \leq p < 2$ or ${\cal B}_0(2,g) > \frac{1}{2 n \pi e} \max_M R_g$, where $R_g$ stands for the scalar curvature of $g$

\item[(VI)] (compactness) ${\cal E}_0(p,g)$ is compact in the $C^1$ topology if either $1 \leq p < 2$ or ${\cal B}_0(2,g) > \frac{1}{2 n \pi e} \max_M R_g$

\end{itemize}

As claimed earlier, local-to-global kind arguments do not seem to be adequate to prove Conjecture I. Instead, using a trick of \cite{Ba}, one concludes that ${\cal A}_0(p,g)$ is well defined for any $1 \leq p < n$. Indeed, given a smooth compact Riemannian manifold $(M,g)$ of dimension $n \geq 2$ and parameters $1 \leq p < n$ and $p \leq q \leq p^*:=\frac{np}{n-p}$, by H\"{o}lder's inequality, one has

\[
||u||_{L^q(M)} \leq ||u||_{L^p(M)}^\alpha ||u||_{L^{p^*}(M)}^{1 - \alpha}\; ,
\]

\n for all $u \in L^{p^*}(M)$, where $\alpha = \frac{np - nq + pq}{pq}$. Taking logarithm of both sides, one gets

\[
\log \left( \frac{||u||_{L^q(M)}}{||u||_{L^p(M)}} \right) + (\alpha - 1) \log \left( \frac{||u||_{L^p(M)}}{||u||_{L^{p^*}(M)}} \right) \leq 0\; .
\]

\n Since this inequality trivializes to an equality when $q = p$, we may differentiate it with respect to $q$ at $q = p$ and it immediately follows that

\[
\int_M |u|^p \log |u|^p \; dv_g \leq \frac{n}{p} \log \left( \int_M |u|^{p^*} \; dv_g\right)^{\frac{p}{p^*}}
\]

\n for all $u \in L^{p^*}(M)$ with $||u||_{L^p(M)} = 1$. So, the continuity of the embedding $H^{1,p}(M) \hookrightarrow L^{p^*}(M)$ yields our claim. Moreover, one easily deduces that

\[
{\cal A}_0(p,g) \leq A_0(p):= \frac{1}{n} \left( \frac{p - 1}{n-p} \right)^{p - 1} \pi^{-\frac{p}{2}}
\left( \frac{\Gamma(n)\Gamma(\frac{n}{2} + 1)}{\Gamma(\frac{n}{p})\Gamma(\frac{n(p - 1)}{p} + 1)} \right)^{\frac{p}{n}} \; ,
\]

\n where $A_0(p)$ is the value of the first best Sobolev constant corresponding to the Riemannian $L^p$-Sobolev inequality.

The main focuses of this work are the first best constant and validity conjectures (or Conjectures II and III).\\

\n 1.3. {\bf Main theorem and proof program} \\

Our main contributions are gathered in the following result:

\begin{teor} \label{MT}
Let $(M,g)$ be a smooth compact Riemannian manifold of dimension $n \geq 2$. For any $1 < p \leq 2$ and $p < n$, we have:

\begin{itemize}

\item[(a)] ${\cal A}_0(p,g) = {\cal A}_0(p)$

\item[(b)] there exists a constant ${\cal B} \in \R$ such that ($L({\cal A}_0(p,g),{\cal B})$) holds for all $u \in H^{1,p}(M)$ with $\int_M |u|^p dv_g = 1$

\end{itemize}

\end{teor}

\n In particular, Conjectures II and III are true or any dimensions $n \geq 3$ and parameters $1 < p \leq 2$.

When we apply the idea of getting entropy inequalities as a limit case of Gagliardo-Nirenberg inequalities to the Riemannian context, the situation changes drastically once extremal functions are usually unknown. In order to solve Conjecture III, Ceccon and Montenegro have examined the validity of sharp Riemannian $L^p$-Gagliardo-Nirenberg inequalities. The proof program of Theorem \ref{MT} relies on these inequalities for which we now provide an overview.

Let $1 \leq p < n$ and $1 \leq q < r < p^*$. An interpolation inequality and the continuity of the embedding $H^{1,p}(M) \hookrightarrow L^{p^*}(M)$ produce constants $A, B \in \R$
such that, for any $u \in H^{1,p}(M)$,

\begin{gather}\label{AB1}
\left( \int_M |u|^r\; dv_g \right)^{\frac{p}{r \theta}} \leq \left( A \int_M |\nabla_g u|^p\; dv_g + B \int_M |u|^p\; dv_g \right) \left(
\int_M |u|^q\; dv_g \right)^{\frac{p(1 - \theta)}{q \theta}}\; ,
\tag{$I_{p,q,r}(A,B)$}
\end{gather}

\n where $\theta = \frac{np(r - q)}{r(q(p - n) + np)} \in (0,1)$.

A few basic notations and definitions associated to (\ref{AB1}) are as follows:

The {\bf first Riemannian $L^p$-Gagliardo-Nirenberg best constant} is defined by

\[
A_0(p,q,r,g) := \inf \{ A \in \R:\; \mbox{ there exists} \hspace{0,18cm} {\cal B} \in \R \hspace{0,18cm} \mbox{such that (\ref{AB})} \hspace{0,18cm} \mbox{holds for all} \hspace{0,18cm} u \in H^{1,p}(M)\}\; .
\]

The {\bf first sharp Riemannain $L^p$-Gagliardo-Nirenberg inequality} states that, there exists a constant $B \in \R$ such
that, for any $u \in H^{1,p}(M)$,

\[
\left( \int_M |u|^r\; dv_g \right)^{\frac{p}{r \theta}} \leq \left(
A_0(p,q,r,g) \int_M |\nabla_g u|^p\; dv_g + B \int_M |u|^p\; dv_g \right)
\left( \int_M |u|^q\; dv_g \right)^{\frac{p(1 - \theta)}{q \theta}}\; .
\]

Whenever this inequality is true, we can define the {\bf second Riemannian $L^p$-Gagliardo-Nirenberg best constant} as

\[
B_0(p,q,r,g) := \inf \{ B \in \R:\; (I_{p,q,r}(A_0(p,q,r,g),B)) \hspace{0,18cm} \mbox{holds for all} \hspace{0,18cm} u \in H^{1,p}(M)\}\; .
\]

\n In this case, the {\bf second sharp Riemannain $L^p$-Gagliardo-Nirenberg inequality} holds on $H^{1,p}(M)$, namely

\[
\left( \int_M |u|^r\; dv_g \right)^{\frac{p}{r \theta}} \leq \left( A_0(p,q,r,g) \int_M |\nabla_g u|^p\; dv_g + B_0(p,q,r,g) \int_M |u|^p\; dv_g \right)
\left( \int_M |u|^q\; dv_g \right)^{\frac{p(1 - \theta)}{q \theta}}\; .
\]

Rearranging the above inequality and taking logarithm of both sides, one obtains

\begin{equation} \label{AE}
\frac{p}{\theta} \log{\frac{||u||_{L^r(M)}}{||u||_{L^q(M)}}} \leq \log{\frac{A_0(p,q,r,g) \int_M |\nabla_g u|^p\; dv_g + B_0(p,q,r,g) \int_M |u|^p\; dv_g}
{\left( \int_M |u|^q\; dv_g \right)^{\frac{p}{q}}} }\, .
\end{equation}

\n The general program of proof consists in finding a suitable class of parameters $q$ and $r$ converging to $p$ such that

\begin{itemize}

\item[($\alpha.1$)] $A_0(p,q,r,g)$ converges to ${\cal A}_0(p,g)$ as $q, r \rightarrow p$

\item[($\alpha.2$)] $B_0(p,q,r,g)$ is bounded for $q$ and $r$ close to $p$

\item[($\alpha.3$)] the inequality (\ref{AE}) yields ($L({\cal A}_0(p,g),{\cal B})$) as $q, r \rightarrow p$

\end{itemize}

The success of this plan relies first on knowing, given $1 < p \leq 2$, if the sharp inequality $(I_{p,q,r}(A_0(p,q,r,g),B))$ is valid for some constant $B$. This question was answered by Ceccon and Montenegro in two works \cite{CM1} and \cite{CM2}. In the first one, in 2008, we prove the validity for $p < r$ and in the second one, in 2013, for $p \geq r$. In practice, the above-described proof strategy is hard due mainly to the difficulty of the assertions ($\alpha.1$) and ($\alpha.2$). In the most cases, the values of $A_0(p,q,r,g)$ and ${\cal A}_0(p,g)$ are unknown. However, one knows that $A_0(p,q,r,g) = A_0(p,q,r)$ whenever $p \leq r$, where $A_0(p,q,r)$ is the best constant to the Euclidean Gagliardo-Nirenberg inequality

\[
\left( \int_{\R^n} |u|^r\; dx \right)^{\frac{p}{r \theta}} \leq \left(
A_0(p,q,r) \int_{\R^n} |\nabla u|^p\; dx\right)
\left( \int_{\R^n} |u|^q\; dx \right)^{\frac{p(1 - \theta)}{q \theta}}\; .
\]

\n This last claim follows directly from a usual unity partition argument.

The major motivation of \cite{CM1} was proving the validity conjecture by choosing parameters introduced by Del Pino and Dolbeault in \cite{DPDo}, namely

\[
s > p \ \ {\rm and}\ \ r = \frac{p(s-1)}{p-1}\, .
\]

\n Indeed, in occasion, we knew that $A_0(p,s,r,g) = A_0(p,s,r)$, since $p < r$, and, as shown in \cite{DPDo}, that $A_0(p,s,r)$ converged to ${\cal A}_0(p)$ as $s \rightarrow p^+$. Moreover, it was possible to prove that ${\cal A}_0(p,g) \geq {\cal A}_0(p)$ and that ($\alpha.3$) followed exactly as in \cite{DPDo}. There remained then the point ($\alpha.2$). However, our strategy was stopped because we were not able to prove the bound of $B_0(p,s,r,g)$ for the above family of parameters $p, s$ and $r$. In other words, we were not able to perform Sections 4 and 5 of this work for such parameters. Yet regarding this family, Brouttelande in \cite{Bro} provided a short proof of ($\alpha.2$) for $p = 2$. His statement assumed that $B_0(p,q,r,g)$ is monotone with respect to $q < 2$, however we were not able to understand why this fact is a direct consequence from his remarks. The author also claims directly that ${\cal A}_0(2,g) = {\cal A}_0(2)$, that is, Conjecture II is always valid for $p = 2$.

The proof of Conjectures II and III were recently resumed after Ceccon and Montenegro \cite{CM2} answered positively that $(I_{p,q,p}(A_0(p,q,p,g),B))$ holds for some constant $B$ whenever $1 < p \leq 2$ and $p < n$. In short, the key ingredients involved in the proof are:

\begin{itemize}

\item[($\beta.1$)] $A_0(p,q,p)$ converges to ${\cal A}_0(p)$ as $q \rightarrow p^-$

\item[($\beta.2$)] $B_0(p,q,p,g)$ is bounded for $q$ close to $p$

\item[($\beta.3$)] ${\cal A}_0(p,g) \geq {\cal A}_0(p)$

\end{itemize}

\n The proof of the assertion ($\beta.1$) uses, among some ideas, a monotonicity result with respect to parameters for the best Gagliardo-Nirenberg constants $A_0(p,q,r)$ of \cite{Ba} and the fact, proved in \cite{DPDo},  that $A_0(p,s,\frac{p(s-1)}{p-1})$ converges to ${\cal A}_0(p)$ as $s \rightarrow p^+$. When $r = p$, we emphasize that extremals for $A_0(p,q,r)$ are unknown. Unlike the first attempt in showing that $B_0(p,s,r,g)$ is bounded for Del Pino and Dolbeault's parameters $s$ and $r$ close to $p$, we here establish the point ($\beta.2$). The proof is done by contradiction and requires concentration and blow-up analysis which it is more delicate than those ones arisen in the Sobolev and Gagliardo-Nirenberg contexts. The proof of the assertion ($\beta.3$) is based on estimates of Gaussian bubbles and Cartan's expansion of metrics around a point on $M$.

Joining ($\beta.1$), ($\beta.2$) and ($\alpha.3$), one discovers that the inequality ($L({\cal A}_0(p),{\cal B})$) is true for some constant ${\cal B} \in \R$. In particular, we have ${\cal A}_0(p,g) \leq {\cal A}_0(p)$, so that, by ($\beta.3$), one concludes the proof of Theorem \ref{MT}.

Section 2 is devoted to the proofs of ($\beta.1$) and ($\beta.3$) and Sections 3, 4 and 5 are devoted to the proof of ($\beta.2$). Thanks to ($\beta.1$), the proof of ($\alpha.3$) is direct.\\

\n 1.4. {\bf Open questions}\\

Assuming Conjecture III to be true for $p=2$, by using Theorem \ref{Ba} of Bakry, Brouttelande \cite{Bro1} showed that, for any $n \geq 2$,

\begin{equation} \label{LB}
{\cal B}_0(2,g) \geq \frac{1}{2 n \pi e} \max_M R_g
\end{equation}

\n and, moreover, that Conjecture V is true in the case that the above inequality is strict. That's all so far known. In particular, it follows completely open that:

\begin{itemize}
\item Conjecture I for any $p \geq n$

\item Conjecture II for $p = 1,  n \geq 2$, $p = 2, n = 2$ and $p > 2, n \geq 2$

\item Conjecture III for $p = 1,  n \geq 2$ and $p = 2, n = 2$

\item Conjectures IV

\item Conjecture V for $1 \leq p < 2,  n \geq 2$

\item Conjecture VI
\end{itemize}

\section{Miscellaneous on best constants}

In this section, we carry all notations of Sobolev type spaces and best constants introduced in the previous section.\\

\n {\bf 2.1 Asymptotic behavior of $A_0(p,q,p)$}. In this subsection, we prove the following result:

\begin{propo}\label{P.1}
For $1 < p < n$ and $1 \leq q < p$, we have

\[
\lim \limits_{q \rightarrow p^-} A_0(p,q,p) = {\cal A}_0(p)\; .
\]

\end{propo}

\n {\bf Proof of Proposition \ref{P.1}.} Inspired on some ideas of \cite{Ba}, we can see that the best constant $A_0(p,q,r)$ is monotone on $q$ and $r$ for each fixed $1 \leq p < n$. Namely, given parameters $q_1 \leq q_2$ and $r_1 \leq r_2$ respecting the range $1 \leq q_i < r_i < p^*$, we have

\begin{equation}\label{firstconst}
A_0(p,q_1,r_1) \leq A_0(p,q_2,r_2)\; .
\end{equation}

\n In fact, it suffices to check that

\[
\left(\int_{\R^n} |u|^{r_1}\; dx \right)^{\frac{p}{r_1 \theta_1}} \leq A_0(p,q_2,r_2) \left( \int_{\R^n} |\nabla u|^p\; dx \right)
\left( \int_{\R^n} |u|^{q_1}\; dx \right)^{\frac{p(1 - \theta_1)}{\theta_1 q_1}}\; ,
\]

\n where $\theta_1 = \frac{np(r_1 - q_1)}{r_1(q_1(p - n) + np)}$.

\n On one hand, a usual interpolation inequality yields

\[
\left( \int_{\R^n} |u|^{r_1}\; dx \right)^{\frac{1}{r_1}} \leq \left(\int_{\R^n} |u|^{q_1}\; dx \right)^{\frac{\lambda}{q_1}}
\left(\int_{\R^n} |u|^{r_2}\; dx \right)^{\frac{1 - \lambda}{r_2}}\; ,
\]

\n where $\frac{1}{r_1} = \frac{\lambda}{q_1} + \frac{1 - \lambda}{r_2}$.

\n In a similar way, we have

\[
\left( \int_{\R^n} |u|^{q_2}\; dx \right)^{\frac{1}{q_2}} \leq \left(\int_{\R^n} |u|^{q_1}\; dx \right)^{\frac{\mu}{q_1}}
\left(\int_{\R^n} |u|^{r_2}\; dx \right)^{\frac{1 - \mu}{r_2}}\; ,
\]

\n where $\frac{1}{q_2} = \frac{\mu}{q_1} + \frac{1 - \mu}{r_2}$.

\n Plugging these two inequalities in

\[
\left(\int_{\R^n} |u|^{r_2}\; dx \right)^{\frac{p}{r_2 \theta_2}} \leq A_0(p,q_2,r_2) \left( \int_{\R^n} |\nabla u|^p\; dx \right)
\left( \int_{\R^n} |u|^{q_2}\; dx \right)^{\frac{p(1 - \theta_2)}{q_2 \theta_2}} \; ,
\]

\n where $\theta_2 = \frac{np(r_2 - q_2)}{r_2(q_2(p - n) + np)}$, one arrives at

\[
\left(\int_{\R^n} |u|^{r_1}\; dx \right)^{\frac{p}{r_1(1 - \lambda)}(1 + \mu\frac{1 - \theta_2}{\theta_2})} \leq A_0(p,q_2,r_2)
\left(\int_{\R^n} |\nabla u|^p\; dx \right) \left( \int_{\R^n} |u|^{q_1}\; dx\right)^{\frac{p}{q_1} (\frac{\mu}{1 - \lambda}
\frac{1 - \theta_2}{\theta_2} + \frac{\lambda}{1 - \lambda})}.
\]

\n On the other hand, from the definition of $\lambda$ and $\mu$ and straightforward computations, we derive

\[
\frac{p}{r_1(1 - \lambda)} \left(1 + \mu \frac{1 -
\theta_2}{\theta_2} \right) = \frac{p}{r_1 \theta_1}
\]

\n and

\[
\frac{p}{q_1}\left( \frac{\mu}{1 - \lambda} \frac{1 -
\theta_2}{\theta_2} + \frac{\lambda}{1 - \lambda}\right) =
\frac{p}{q_1}\frac{1 - \theta_1}{\theta_1} \; ,
\]

\n so that the assertion (\ref{firstconst}) follows.

\n For fixed $1 < p < n$, choose parameters $1 \leq q < p$ and $p < s <\frac{p (n - 1)}{n - p}$. From the previous conclusion, we have

\[
A_0(p,q,p) \leq A_0(p,s, \frac{p(s - 1)}{p - 1}) \; .
\]

\n According to \cite{DPDo}, the right-hand side converges to ${\cal A}_0(p)$ as $s \rightarrow p^+$, so that

\[
\limsup \limits_{q \rightarrow p^-} A_0(p,q,p) \leq {\cal A}_0(p) \; .
\]

\n By an adaptation of an argument of \cite{Ba}, we next show that

\[
\liminf \limits_{q \rightarrow p^-} A_0(p,q,p) \geq {\cal A}_0(p) \; .
\]

\n In fact, consider the sharp Euclidean Gagliardo-Nirenberg inequality

\[
\left( \int_{\R^n} |u|^p \; dx \right)^{\frac{1}{\theta_q}} \leq A_0(p,q,p) \left( \int_{\R^n} |\nabla u|^p \; dx \right) \left(
\int_{\R^n} |u|^q \; dx \right)^{\frac{p(1 - \theta_q)}{q \theta_q}}\; ,
\]

\n where $\theta_q = \frac{n(p - q)}{np + pq - nq} \in (0,1)$. Taking logarithm of both sides, one has

\[
\frac{1}{\theta_q} \log \left(\frac{||u||_p}{||u||_q} \right) \leq \log \left(A_0(p,q,p)\frac{||\nabla u||_p^p}{||u||_q^p}\right)^{\frac{1}{p}}\; .
\]

\n Using the definition of $\theta_q$ and taking the limit on $q$, one has

\[
\frac{p^2}{n} \lim_{q \rightarrow p^-} \frac{1}{p - q} \log \left(\frac{||u||_p}{||u||_q} \right)\leq \log \left(\liminf \limits_{q \rightarrow p^-} A_0(p,q,p)
\frac{||\nabla u||_p^p}{||u||_p^p}\right)^{\frac{1}{p}}\; .
\]

\n We now compute the left-hand side limit. First write

\[
\log \left(\frac{||u||_p}{||u||_q} \right) = \frac{1}{p} \log(||u||_p^p) - \frac{1}{q} \log(||u||_q^q) = \frac{q - p}{p}
\log(||u||_q) + \frac{1}{p} \left( \log(||u||_p^p) - \log(||u||_q^q) \right)\; .
\]

\n Applying two times the mean value theorem, one gets

\[
\lim_{q \rightarrow p^-} \frac{1}{p - q} \log \left(\frac{||u||_p}{||u||_q} \right) = \frac{1}{p} \int_{\R^n}
\frac{|u|^p}{||u||_p^p} \log \left(\frac{|u|}{||u||_p}\right)\; dx\; .
\]

\n Therefore,

\[
\int_{\R^n} \frac{|u|^p}{||u||_p^p} \log \left( \frac{|u|^p}{||u||_p^p} \right)\; dx \leq \frac{n}{p} \log \left( \liminf \limits_{q \rightarrow p^-} A_0(p,q,p) \frac{||\nabla
u||_p^p}{||u||_p^p}\right)\; ,
\]

\n or equivalently,

\[
\int_{\R^n}|u|^p \log (|u|^p)\; dx \leq \frac{n}{p} \log \left( \liminf \limits_{q \rightarrow p^-} A_0(p,q,p)
\int_{\R^n} |\nabla u|^p\; dx \right)
\]

\n for all $u \in C^\infty_0(\R^n)$ with $||u||_p = 1$. But just this inequality implies the desired assertion.\bl\\

\n {\bf 2.2 A lower bound for ${\cal A}_0(p,g)$}. As mentioned in the introduction, we need the following result:

\begin{propo}\label{P.2}
For each $n \geq 2$ and $1 < p < n$, we have

\[
{\cal A}_0(p,g) \geq {\cal A}_0(p)\; .
\]

\end{propo}

\n {\bf Proof of Proposition \ref{P.2}.}  Let ${\cal A}$ and ${\cal B}$ be constants such that $L({\cal A},{\cal B})$ is valid. It suffices to show that ${\cal A} \geq {\cal A}_0(p)$. We proceed by contradiction. Assume by contradiction that ${\cal A} < {\cal A}_0(p)$. One knows that $L({\cal A},{\cal B})$ is equivalent to

\begin{equation} \label{EF1}
\frac{1}{||u||_p^p} \int_M |u|^p \log(|u|^p)\; dv_g + (\frac{n}{p} - 1) \log (||u||_p^p) \leq \frac{n}{p} \log \left(  {\cal A} \int_M |\nabla
u|_g^p\; dv_g + {\cal B} \int_M |u|^p\; dv_g \right)
\end{equation}

\n for all $u \in C^{\infty}(M)$.

\n Let us fix a point $x_0 \in M$ and an extremal function $u_0(x) = a e^{-b|x|^{\frac{p}{p-1}}} \in W^{1,p}(\R^n)$ for the sharp Euclidean entropy inequality

\[
\int_{\R^n} |u|^p \log(|u|^p)\; dx \leq \frac{n}{p} \log \left( {\cal A}_0(p) \int_M |\nabla u|^p\; dx \right)\; ,
\]

\n where $a$ and $b$ are positive constants such that $||u_0||_p = 1$, see [\cite{DPDo} or \cite{G}]. Choose now a geodesic ball $B(x_0,\delta) \subset M$ and a radial cutoff function $\eta \in C^\infty(B(0,\delta))$ satisfying $\eta = 1$ in $B(0,\frac{\delta}{2})$, $\eta = 0$ outside $B(0,\delta),\ 0 \leq \eta \leq 1$ in $B(0,\delta)$. For $\varepsilon > 0$ and $x \in B(0,\delta)$, set

\[
u_\varepsilon(exp_{x_0}(x)) = \eta(x) \varepsilon^{-\frac{n}{p}} u_0(\frac{x}{\varepsilon})\; .
\]

\n The asymptotic behavior of some integrals computed for $u_\varepsilon$ with $\varepsilon > 0$ small enough are now presented. Denote

\[
I_1 = \int_{\R^n} u_0^p \log(u_0^p)\; dx, \ \ I_2 = \int_{\R^n} |\nabla u_0|^p\; dx\; ,
\]

\[
J_1 = \int_{\R^n} u_0^p |x|^2\; dx, \ \ J_2 = \int_{\R^n} |\nabla u_0|^p |x|^2\; dx,\ \ J_3 = \int_{\R^n} u_0^p \log(u_0^p) |x|^2\; dx\; .
\]

\n Using the expansion of volume element in geodesic coordinates

\[
\sqrt{det g} = 1 - \frac{1}{6} \sum_{i,j = 1}^{n} Ric_{ij}(x_0) x_i x_j + O(r^3)\; ,
\]

\n where $Ric_{ij}$ denotes the components of the Ricci tensor in these coordinates and $r=|z|$, one easily checks that

\begin{equation} \label{1}
\int_{M} u_\varepsilon^p\; dv_g = 1 - \frac{R_g(x_0)}{6n} J_1 \varepsilon^2 + O(\varepsilon^4)\; ,
\end{equation}

\begin{equation} \label{2}
\int_{M} u_\varepsilon^p \log(u_\varepsilon^p)\; dv_g = I_1 - n \log \varepsilon - \frac{R_g(x_0)}{6n} J_3 \varepsilon^2 + \frac{R_g(x_0)}{6} J_1 \varepsilon^2 \log \varepsilon  + o(\varepsilon^4 \log \varepsilon)
\end{equation}

\n and

\begin{equation} \label{3}
\int_{M} |\nabla u_\varepsilon|^p\; dv_g = \varepsilon^{-p} \left(I_2 - \frac{R_g(x_0)}{6n} J_2 \varepsilon^2 + O(\varepsilon^4) \right)\; .
\end{equation}

\n Plugging $u_\varepsilon$ in (\ref{EF1}) and evoking the asymptotic behaviors (\ref{1}), (\ref{2}) and (\ref{3}), one obtains

\[
\frac{I_1 - n \log \varepsilon - \frac{R_g(x_0)}{6n} J_3 \varepsilon^2 + \frac{R_g(x_0)}{6} J_1 \varepsilon^2 \log \varepsilon  + o(\varepsilon^4 \log \varepsilon)}{1 - \frac{R_g(x_0)}{6n} J_1 \varepsilon^2 + O(\varepsilon^4)} + (\frac{n}{p} - 1) \log \left(1 - \frac{R_g(x_0)}{6n} J_1 \varepsilon^2 + O(\varepsilon^4)\right)
\]

\[
\leq -n \log \varepsilon + \frac{n}{p} \log\left( {\cal A} I_2 - \frac{R_g(x_0)}{6n} {\cal A} J_2 \varepsilon^2 + {\cal B} \varepsilon^p +
O(\varepsilon^q)\right)
\]

\n for $\varepsilon > 0$ small enough, where $q = \min\{4, p +2\}$.

\n So, Taylor`s expansion guarantees that

\begin{equation}\label{taylor}
I_1 - n \log \varepsilon - \frac{R_g(x_0)}{6n} J_3 \varepsilon^2 +
\frac{R_g(x_0)}{6} J_1 \varepsilon^2 \log \varepsilon +
\frac{R_g(x_0)}{6n} I_1 J_1 \varepsilon^2 - \frac{R_g(x_0)}{6} J_1
\varepsilon^2 \log \varepsilon
\end{equation}

\[
- (\frac{n}{p} - 1) \frac{R_g(x_0)}{6n} J_1 \varepsilon^2 +
O(\varepsilon^4 \log \varepsilon) \leq - n \log \varepsilon +
\frac{n}{p} \log ({\cal A} \; I_2) - \frac{n}{p} \frac{R_g(x_0)}{6n}
\frac{J_2}{I_2} \varepsilon^2 + \frac{n}{p} \frac{{\cal B}}{{\cal A} \; I_2}
\varepsilon^p + O(\varepsilon^q)\; .
\]

\n Thanks to the assumption by contradiction ${\cal A} < {\cal A}_0(p)$, after a suitable simplification, one arrives at the following contradiction

\[
0 < I_1 - \frac{n}{p} \log({\cal A} \; I_2) \leq - \frac{R_g(x_0)}{6n} \left( I_1 J_1 + \frac{n}{p} \frac{J_2}{I_2} - (\frac{n}{p} - 1) J_1 -
J_3\right) \varepsilon^2 + O(\varepsilon^p)
\]

\n for $\varepsilon > 0$ small enough.\bl\\

\section{The bound of $B_0(p,q,p,g)$}

The remaining of this work is devoted to the proof of the following theorem:

\begin{teor}\label{lB}
Let $(M,g)$ be a smooth compact Riemannian manifold of dimension $n \geq 2$. For each fixed $1 < p \leq 2$ and $p < n$, the best constant $B_0(p,q,p,g)$ is
bounded for any $q < p$ close to $p$.
\end{teor}

In this section, we present the PDEs setting after assuming the above assertion fails. Namely, given $1 < p \leq 2$ and $p < n$, suppose by contradiction that the conclusion of Theorem \ref{lB} is false. Equivalently, by Proposition \ref{P.1}, there exists a sequence $(q) \subset (1,p)$ such that

\[
\lim \limits_{q \rightarrow p^-} C_q = + \infty\; ,
\]

\n where

\[
C_q = \frac{B_0(p,q,p,g) - (p-q)}{A_0(p,q,p)}\; .
\]

\n Since $C_q < B_0(p,q,p,g) A_0(p,q,p)^{-1}$, we have

\begin{equation}\label{ldf}
\nu_q = \inf_{u \in {\cal H}} J_q (u) < A_0(p,q,p)^{-1}\; ,
\end{equation}

\n where ${\cal H} = \{ u \in H^{1,p}(M):\; ||u||_{L^p(M)} = 1 \}$ and

\[
J_q(u) = \left( \int_M |\nabla_g u|^p\; dv_g + C_q \int_M |u|^p\; dv_g \right) \left( \int_M |u|^q\; dv_g \right)^{\frac{p(1 -
\theta_q)}{q \theta_q}}\; .
\]

\n Note that $J_q$ is of class $C^1$, so that standard variational methods involving constrained minimization easily produce a minimizer $u_q \in {\cal H}$ of $J_q$, that is

\begin{equation}\label{l3iha}
J_q(u_q) = \nu_q = \inf_{u \in {\cal H}} J_q(u)\; .
\end{equation}

\n One may assume $u_q \geq 0$, since $|\nabla_g |u_q|| = |\nabla_g u_q|$ a.e.. Moreover, as can be easily checked, $u_q$ satisfies

\begin{equation}\label{l3ep}
A_q \Delta_{p,g} u_q + A_q C_q u_q^{p - 1} + \frac{1 -
\theta_q}{\theta_q} B_q u_q^{q - 1} = \frac{\nu_q}{\theta_q} u_q^{p - 1}
\hspace{0,2cm} \mbox{on} \hspace{0,2cm} M\; ,
\end{equation}

\n where $\Delta_{p,g} = -{\rm div}_g(|\nabla_g|^{p-2} \nabla_g)$ is the $p$-Laplace operator of $g$,

\[
A_q = \left(\int_M u_q^q\; dv_g \right)^{\frac{p( 1 - \theta_q)}{q \theta_q}}\; ,
\]

\[
B_q = \left(\int_M |\nabla_g u_q|^p\; dv_g + C_q \right) \left(\int_M u_q ^q\; dv_g \right)^{\frac{p(1 - \theta_q)}{q \theta_q} - 1}\; .
\]

\n In particular, we have the relation

\begin{equation} \label{nu}
B_q \int_M u_q^q\; dv_g =\nu_q\; .
\end{equation}

\n By the Tolksdorf's regularity theory \cite{To}, it follows that $u_q$ is of class $C^1$.

\n We now assert that

\begin{equation} \label{asym}
\lim \limits_{q \rightarrow p^-} \nu_q = {\cal A}_0(p)^{-1}\; .
\end{equation}

\n Indeed, a combination between the validity of $I_{p,q,p}(A_0(p,q,p),B_0(p,q,p,g))$ and the definition of $C_q$ gives

\[
A_0(p,q,p)^{-1} \leq \left(\int_M |\nabla_g u_q|^p\; dv_g + C_q \right)\left(\int_M u_q^q\; dv_g \right)^{\frac{p(1 - \theta_q)}{q \theta_q}} + \frac{(p - q)}{A_0(p,q,p)} \; A_q  = \nu_q + A_q \;\frac{p - q}{A_0(p,q,p)}\; .
\]

\n By Proposition \ref{P.1}, $A_0(p,q,p)$ remains away from zero when $q$ is close to $p$. So,

\[
\liminf \limits_{q \rightarrow p^-} \nu_q \geq {\cal A}_0(p)^{-1}
\]

\n and the conclusion follows from the inequality $\nu_q < A(p,q,p)^{-1}$ for all $1 \leq q < p$.\bl \\

\section{Concentration analysis of $(u_q)$}

\n Let $x_q \in M$ be a maximum point of $u_q$, that is

\begin{equation}\label{l3ix}
u_q(x_q) = ||u_q||_{L^\infty(M)}\; .
\end{equation}

\n Our aim here is to establish the following concentration property satisfied by $(u_q)$:

\begin{lema}\label{CA}
We have

\[
\lim \limits_{\sigma\rightarrow +\infty} \lim \limits_{q \rightarrow p^-} \int_{B(x_q,\sigma A_q^{\frac{1}{p}})} u_q^p\; dv_g = 1 \; .
\]
\end{lema}

\n {\bf Proof of Lemma \ref{CA}.} By (\ref{ldf}), one has $A_q C_q < A_0(p,q,p)^{-1}$, so that, by Proposition \ref{P.1}, $A_q \rightarrow 0$ as $q \rightarrow p^-$.

\n Let $\sigma > 0$. For each $x \in B(0, \sigma)$, define

\[
\begin{array}{l}
h_q(x) = g(\exp_{x_q} (A_q^{\frac{1}{p}} x))\; , \vspace{0,3cm}\\
\varphi_q(x) = A_q ^{\frac{n}{p^2}}
u_q(\exp_{x_q}(A_q^{\frac{1}{p}} x))\; .
\end{array}
\]

\n By (\ref{l3ep}) and (\ref{nu}), one easily deduces that

\[
\Delta_{p,h_q} \varphi_q + A_q C_q \varphi_q^{p - 1} + \frac{1 - \theta_q}{\theta_q} \nu_q \varphi_q^{q - 1} = \frac{\nu_q}{\theta_q} \varphi_q^{p - 1} \hspace{0,2cm} \mbox{on} \hspace{0,2cm}
B(0,\sigma) \; .
\]

\n Using the mean value theorem and the value of $\theta_q$, one gets

\begin{equation}\label{el1}
\Delta_{p,h_q} \varphi_q + A_q C_q \varphi_q^{p - 1} = \nu_q ( \varphi_q^{q - 1} + \frac{pq - nq + np}{n} \varphi_q^{\rho_q} \log (\varphi_q)) \hspace{0,2cm} \mbox{on} \hspace{0,2cm} B(0,\sigma)
\end{equation}

\n for some $\rho_q \in (q - 1,p - 1)$.

\n Consider $\varepsilon > 0$ fixed such that $p + \varepsilon < \frac{np}{n - p}$. Since $\varphi_q^{\rho_q} \log (\varphi_q^\varepsilon) \leq \varphi_q^{p - 1 + \varepsilon}$, we then get

\[
\Delta_{p,h_q} \varphi_q + A_q C_q \varphi_q^{p - 1} \leq \nu_q ( \varphi_q^{q - 1} + \frac{pq - nq + np}{n \varepsilon}
\varphi_q^{\varepsilon + p - 1}) \hspace{0,2cm} \mbox{on} \hspace{0,2cm} B(0,\sigma)\; .
\]

\n Since all coefficients of this equation are bounded, the Moser's iterative scheme (see \cite{Se}) produces

\[
(A_q^{\frac{n}{p^2}} ||u_q||_{L^{\infty}(M)})^p = \sup_{B(0,\sigma)} \varphi_q^p \leq c_\sigma \int_{B(0,2\sigma)} \varphi_q^p\; dh_q =
c_\sigma \int_{B(x_q,2 \sigma A_q^{\frac{1}{p}})} u_q^p\; dv_g \leq c_\sigma
\]

\n for all $q$ close to $p$, where $c_\sigma$ is a constant independent of $q$.

\n On the other hand,

\[
1 = ||u_q||_{L^p(M)}^p \leq ||u_q||_{L^{\infty}(M)}^{p-q} \, ||u_q||_{L^q(M)}^q = \left(||u_q||_{L^{\infty}(M)} \, A_q^{\frac{n}{p^2}}\right)^{p - q} \; ,
\]

\n so that

\begin{equation}\label{lep}
1 \leq ||\varphi_q||_{L^{\infty}(B(0,\sigma))} \leq c_\sigma^{1/p}
\end{equation}

\n for all $q$ close to $p$.

\n By (\ref{ldf}) and Proposition \ref{P.1}, up to a subsequence, we have

\[
{\cal A}_0(p)^{-1} \geq \lim_{q \rightarrow p^-} A_q C_q = C \geq 0\; .
\]

\n Thanks to (\ref{lep}) and Tolksdorf's elliptic theory applied to (\ref{el1}), it follows that $\varphi_q \rightarrow \varphi$ in $C^1_{loc}(\R^n)$ and $\varphi \not \equiv 0$. Letting $q \rightarrow p^-$ in (\ref{el1}) and using (\ref{asym}), one has

\begin{equation} \label{ELim}
\Delta_{p,\xi} \varphi + C \varphi^{p - 1} = {\cal A}_0(p)^{-1} \left( \varphi^{p - 1} + \frac{p}{n} \varphi^{p - 1} \log (\varphi^p)\right) \hspace{0,2cm} \mbox{on} \hspace{0,2cm} \R^n \;,
\end{equation}

\n where $\Delta_{p,\xi}$ stands for the Euclidean $p$-Laplace operator.

\n For each $\sigma > 0$, we have

\begin{equation} \label{whole}
\int_{B(0,\sigma)} \varphi^p\; dx = \lim_{q \rightarrow p^-} \int_{B(0,\sigma)} \varphi_q^p\; dh_q = \lim_{q \rightarrow p^-}\int_{B(x_q,\sigma A_q^{\frac{1}{p}})} u_q^p\; dv_g \leq 1
\end{equation}

\n and, by $(\ref{ldf})$ and Proposition \ref{P.1},

\begin{equation} \label{grad}
\int_{B(0,\sigma)} |\nabla \varphi|^p\; dx = \lim_{q \rightarrow p^-} \int_{B(0,\sigma)} |\nabla _{h_q} \varphi_q|^p\; dh_q = \lim_{q \rightarrow p^-} \left( A_q \int_{B(x_q,\sigma A_q^{\frac{1}{p}})} |\nabla_g u_\alpha|^p\; dv_g \right) \leq {\cal A}_0(p)^{-1}\; .
\end{equation}

\n In particular, one has $\varphi \in W^{1,p}(\R^n)$. Consider then a sequence of nonnegative functions $(\varphi_k) \subset C^\infty_0(\R^n)$ converging to $\varphi$ in $W^{1,p}(\R^n)$. Taking $\varphi_k$ as a test function in (\ref{ELim}), we can write

\[
\frac{n}{p} {\cal A}_0(p) \int_{\R^n} |\nabla \varphi|^{p-2} \nabla \varphi \cdot \nabla \varphi_k\; dx + \frac{n}{p} {\cal A}_0(p) C \int_{\R^n} \varphi^{p-1} \varphi_k\; dx = \int_{\R^n} \varphi^{p-1} \varphi_k \log (\varphi^p)\; dx + \frac{n}{p} \int_{\R^n} \varphi^{p-1} \varphi_k \; dx
\]

\[
= \int_{[\varphi \leq 1]} \varphi^{p-1} \varphi_k \log (\varphi^p)\; dx + \int_{[\varphi \geq 1]} \varphi^{p-1} \varphi_k \log (\varphi^p)\; dx + \frac{n}{p} \int_{\R^n} \varphi^{p-1} \varphi_k \; dx
\]

\n Letting $k \rightarrow + \infty$ and applying Fatou's lemma and the dominated convergence theorem in the right-hand side, one gets

\[
\frac{n}{p} {\cal A}_0(p) \int_{\R^n} |\nabla \varphi|^p\; dx \leq \int_{[\varphi \leq 1]} \varphi^p \log (\varphi^p)\; dx + \int_{[\varphi \geq 1]} \varphi^p \log (\varphi^p)\; dx + \frac{n}{p} \int_{\R^n} \varphi^p\; dx
\]

\[
= \int_{\R^n} \varphi^p \log (\varphi^p)\; dx + \frac{n}{p} \int_{\R^n} \varphi^p\; dx\; .
\]

\n Rewriting this inequality in function of $\psi(x) = \frac{\varphi(x)}{||\varphi||_{L^p(M)}}$, one has

\[
\frac{n}{p} {\cal A}_0(p) \int_{\R^n} |\nabla \psi|^p\; dx \leq \int_{\R^n} \psi^p \log (\psi^p)\; dx + \frac{n}{p} \int_{\R^n} \psi^p\; dx + \log ||\varphi||_{L^p(M)}^p\; .
\]

\n Note that this inequality combined with

\[
\int_{\R^n} \psi^p \log (\psi^p)\; dx  \leq \frac{n}{p} \log \left({\cal A}_0(p) \int_{\R^n} |\nabla \psi|^p\; dx \right)
\]

\n produces

\[
{\cal A}_0(p) \int_{\R^n} |\nabla \psi|^p\; dx \leq \log \left({\cal A}_0(p)\; \int_{\R^n} |\nabla \psi|^p\; dx\right) + 1 + \frac{p}{n} \log ||\varphi||_{L^p(M)}^p\; .
\]

\n Since $\log x \leq x - 1$ for all $x > 0$, it then follows that $\log ||\varphi||_{L^p(M)}^p \geq 0$ or, in other words, $||\varphi||_{L^p(M)} \geq 1$. Thus, by (\ref{whole}), we conclude that $||\varphi||_{L^p(M)} = 1$, so that

\[
\lim_{\sigma \rightarrow \infty} \lim_{q \rightarrow p^-} \int_{B(x_q,\sigma A_q^{\frac{1}{p}})} u_q^p\; dv_g = \int_{\R^n} \varphi^p\; dx = 1\; .
\]
\bl \\

\section{A uniform estimate for $(u_q)$}

As mentioned in the introduction, an essential lemma in the proof of Theorem \ref{MT} is the following

\begin{lema}\label{UE}
For any $\lambda > 0$, there exists a constant $c_\lambda > 0$, independent of $q < p$, such that

\[
d_g(x,x_q)^\lambda u_q(x) \leq c_\lambda A_q^{\frac{\lambda}{p} - \frac{n}{p^2}}
\]

\n for all $x \in M$ and $q$ close enough to $p$, where $d_g$ stands for the distance with respect to the metric $g$.
\end{lema}

\n {\bf Proof of Lemma \ref{UE}.} Suppose by contradiction that the conclusion of the above lemma fails. Then, there exist $\lambda_0 > 0$ and $y_q \in M$ such that $f_{q,\lambda_0}(y_q) \rightarrow + \infty$ as $q \rightarrow p^-$, where

\[
f_{q,\lambda}(x) = d_g(x,x_q)^\lambda u_q(x) A_q^{\frac{n}{p^2} - \frac{\lambda}{p}} \; .
\]

\n Without loss of generality, assume that $f_{q,\lambda_0}(y_q) = ||f_{q,\lambda_0}||_{L^\infty(M)}$. From (\ref{lep}), we have

\[
f_{q,\lambda_0}(y_q) \leq c \frac{u_q(y_q)}{||u_q||_{L^\infty(M)}} d_g(x_q,y_q)^{\lambda_0} A_q^{- \frac{\lambda_0}{p}} \leq c d_g(x_q,y_q)^{\lambda_0} A_q^{- \frac{\lambda_0}{p}}\; ,
\]

\n so that

\begin{equation}\label{inf}
d_g(x_q,y_q) A_q^{- \frac{1}{p}} \rightarrow + \infty\ \ {\rm as}\ \ q \rightarrow p^-\; .
\end{equation}

\n For any fixed $\varepsilon \in (0,1)$ and $\sigma > 0$, we next show that

\begin{equation}\label{3int}
B(y_q,\varepsilon d(x_q,y_q)) \cap B(x_q, \sigma A_q^{\frac{1}{p}}) = \emptyset
\end{equation}

\n for $q$ close enough to $p$. For this, it suffices to check that

\[
d_g(x_q,y_q) \geq \sigma  A_q^{\frac{1}{p}} + \varepsilon d(x_q,y_q) \; .
\]

\n But the above inequality is equivalent to

\[
d_g(x_q, y_q)(1 - \varepsilon) A_q^{- \frac{1}{p}} \geq \sigma \; ,
\]

\n which is clearly satisfied, since $d_g(x_q,y_q) A_q^{- \frac{1}{p}} \rightarrow + \infty$ for $q$ close to $p$ and $1 - \varepsilon > 0$.

\n We claim that exists a constant $c > 0$, independent of $q$, such that

\begin{equation}\label{lim}
u_q(x) \leq c u_q(y_q)
\end{equation}

\n for all $x \in B(y_q, \varepsilon d_g(x_q,y_q))$ and $q$ close to $p$. Indeed,

\[
d_g(x,x_q) \geq d_g(x_q,y_q) - d_g(x,y_q) \geq (1 - \varepsilon) d_g(x_q,y_q)
\]

\n for all $x \in B(y_q, \varepsilon d_g(x_q,y_q))$.

\n Thus,

\[
d_g(x_q,y_q)^{\lambda_0} u_q(y_q) A_q^{\frac{n}{p^2} - \frac{\lambda_0}{p}} = f_{q,\lambda_0}(y_q) \geq f_{q,\lambda_0}(x) = d_g(x,x_q)^{\lambda_0} u_q(x) A_q^{\frac{n}{p^2} -
\frac{\lambda_0}{p}}
\]

\[
\geq (1 - \varepsilon)^{\lambda_0} d_g(x_q,y_q)^{\lambda_0} u_q(x) \; A_q^{\frac{n}{p^2} - \frac{\lambda_0}{p}} \; ,
\]

\n so that

\[
u_q(x) \leq \left(\frac{1}{1 - \varepsilon}\right)^{\lambda_0} u_q(y_q)
\]

\n for all $x \in B(y_q, \varepsilon d_g(x_q,y_q))$ and $q$ close to $p$, as claimed.

\n Now define

\[
\begin{array}{l}
\tilde{h}_q(x) = g(\exp_{y_q} (A_q^{\frac{1}{p}} x))\; , \vspace{0,3cm}\\
\tilde{\varphi}_q(x) = A_q ^{\frac{n}{p^2}} u_q(\exp_{y_q}(A_q^{\frac{1}{p}} x)) \; .
\end{array}
\]

\n By (\ref{l3ep}) and (\ref{nu}), it readily follows that

\[
\Delta_{p,h_q} \tilde{\varphi}_q + A_q C_q \tilde{\varphi}_q^{p - 1} + \frac{1 - \theta_q}{\theta_q} \nu_q \; \tilde{\varphi}_q^{q - 1}
= \frac{\nu_q}{\theta_q} \tilde{\varphi}_q^{p - 1} \hspace{0,2cm} \mbox{on} \hspace{0,2cm} B(0,3) \; .
\]

\n Applying the mean value theorem, one obtains

\begin{equation}\label{leqgn}
\Delta_{p,h_q} \tilde{\varphi}_q + A_q C_q \tilde{\varphi}_q^{p - 1} = \nu_q \left( \tilde{\varphi}_q^{q - 1} + \frac{pq - nq +
np}{n} \log (\tilde{\varphi}_q) \tilde{\varphi}_q^{\rho_q}\right) \hspace{0,2cm} \mbox{on} \hspace{0,2cm} B(0,3) \; ,
\end{equation}

\n where $\rho_q \in (q - 1,p - 1)$.

\n For fixed $\varepsilon > 0$ such that $p + \varepsilon < \frac{np}{n - p}$, we have $\log (\tilde{\varphi}_q^\varepsilon)
\tilde{\varphi}_q^{\rho_q} \leq \tilde{\varphi}_q^{p - 1 + \varepsilon}$. So, the Moser's iterative scheme applied to (\ref{leqgn}) yields

\begin{equation}\label{pf}
\mu_q^{\frac{p}{p - q}} = (A_q^{\frac{n}{p^2}} u_q(y_q))^p \leq \sup_{B(0,1)} \tilde{\varphi}_q^p \leq c \int_{B(0,2)}
\tilde{\varphi}_q^p\; d\tilde{h}_q = c \int_{B(y_q,A_q^{\frac{1}{p}})} u_q^p\; dv_g
\end{equation}

\n for $q$ close to $p$, where

\begin{equation}\label{mu}
\mu_q = u_q(y_q)^{p - q} \int_M u_q^q dv_g\; .
\end{equation}

\n Note that two independent situations can occur: either

\n (i) there exists a subsequence $(q_k)$ of $(q)$ such that $\mu_{q_k} \geq 1 - \theta_{q_k}$, or

\n (ii) $\mu_q < 1 - \theta_q$ for $q$ close enough to $p$.

\n In each case, we derive a contradiction. In fact, if (i) holds, we have

\[
\liminf \limits_{k \rightarrow + \infty} \mu_{q_k}^{\frac{p}{p - q_k}} \geq e^{-\frac np}\; .
\]

\n On the other hand, by (\ref{inf}), one gets

\begin{equation}\label{sub}
B(y_{q_k}, A_{q_k}^{\frac{1}{p}}) \subset B(y_{q_k}, \varepsilon d(x_{q_k},y_{q_k}))
\end{equation}

\n for $k$ large enough. So, by Lemma \ref{CA}, (\ref{3int}) and (\ref{pf}), one arrives at the contradiction

\[
0 < e^{-\frac np} \leq \lim_{k \rightarrow + \infty} \int_{B(y_{q_k},A_{q_k}^{\frac{1}{p}})} u_{q_k}^p\; dv_g = 0 \; .
\]

\n If the assertion (ii) is satisfied, we set

\[
\begin{array}{l}
\tilde{h}_q(x) = g(\exp_{y_q}(A_q^{\frac{1}{p}} x)) \vspace{0,2cm}\\
\psi_q(x) = u_q(y_q)^{-1} u_q(\exp_{y_q}(A_q^{\frac{1}{p}} x))\; .
\end{array}
\]

\n Thanks to (\ref{l3ep}) and (\ref{nu}), we have

\[
\Delta_{p,h_q} \psi_q + A_q C_q \psi_q^{p - 1} + \frac{1 - \theta_q}{\theta_q} \frac{\nu_q}{\mu_q} \psi_q^{q - 1} =
\frac{\nu_q}{\theta_q} \psi_q^{p - 1} \hspace{0,2cm} \mbox{on} \hspace{0,2cm} B(0,3) \; ,
\]

\n where $\mu_q$ is given in (\ref{mu}). Rewriting this equation as

\[
\Delta_{p,h_q} \psi_q + A_q C_q \psi_q^{p - 1} + \frac{\nu_q}{\theta_q}\left(\frac{1 - \theta_q}{\mu_q} -1\right)
\psi_q^{q - 1} = \frac{\nu_q}{\theta_q}\left( \psi_q^{p - 1} - \psi_q^{q - 1}\right) \hspace{0,2cm} \mbox{on} \hspace{0,2cm} B(0,3) \; ,
\]

\n the mean value theorem gives

\begin{equation} \label{norm}
\Delta_{p,h_q} \psi_q + A_q C_q \psi_q^{p - 1} + \frac{\nu_q}{\theta_q}\left(\frac{1 - \theta_q}{\mu_q} -1\right) \psi_q^{q - 1} = \frac{\nu_q (pq - nq + np)}{n} \log (\psi_q)
\psi_q^{\rho_q}  \hspace{0,2cm} \mbox{on} \hspace{0,2cm} B(0,3)
\end{equation}

\n for some $\rho_q \in (q-1,p-1)$.

\n Using the fact that $\frac{1 - \theta_q}{\mu_q} - 1 > 0$, (\ref{inf}) and (\ref{lim}), one easily deduces that $\psi_q \rightharpoondown \psi$ in
$W^{1,p}(B(0,2))$. Let $h \in C^1_0(B(0,2))$ be such that $h \equiv 1$ in $B(0,1)$ and $h \geq 0$. Since, by the Moser's iterative scheme, $\psi \not \equiv 0$, choosing $\psi h^p$ as a test function in (\ref{norm}), one obtains

\[
\limsup_{q \rightarrow p^-} \frac{1}{\theta_q}\left(\frac{1 - \theta_q}{\mu_q} - 1\right) < c\; .
\]

\n Therefore, up to a subsequence, we can write

\[
\lim_{q \rightarrow p^-} \frac{1}{\theta_q}\left(\frac{1 -
\theta_q}{\mu_q} - 1\right) = \gamma \geq 0\; ,
\]

\n so that $\mu_q \rightarrow 1$. Using the definition of $\theta_q$, we then derive

\[
\lim \limits_{q \rightarrow p^-} \mu_q^{\frac{p}{p - q}} = e^{-(1 + \gamma)\frac{n}{p}}\; .
\]

\n Finally, combining Lemma \ref{CA}, (\ref{3int}), (\ref{pf}) and the above limit, we are led to the contradiction

\[
0 <  e^{-(1 + \gamma)\frac{n}{p}} \leq \lim_{q \rightarrow p^-} \int_{B(y_q,A_q^{\frac{1}{p}})} u_q^p\; dv_g = 0\; .
\]

\bl\\

\section{Achieving the final contradiction}

In the sequel, we perform several estimates with the aid of Lemma \ref{UE} in order to establish the desired contradiction in the proof of Theorem \ref{lB}.

From now on, several possibly different positive constants independent of $q$ will be denoted by
$c$ or $c_1$.

Assume, without loss of generality, that the radius of injectivity of $M$ is greater than $2$. Let $\eta \in C^1_0(\R)$ be a cutoff function such that $\eta = 1$ on $[0,1)$, $\eta = 0$ on $[2, \infty)$ and $0 \leq \eta \leq 1$ and define $\eta_q(x) = \eta(d_g(x,x_q))$.

The sharp Euclidean Gagliardo-Nirenberg inequality provides

\[
\left( \int_{B(0,2)} (u_q \eta_q)^p\; dx \right)^{\frac{1}{\theta_q}} \leq A_0(p,q,p) \left( \int_{B(0,2)} |\nabla (u_q \eta_q)|^p\; dx\right) \left(\int_{B(0,2)} (u_q \eta_q)^q\; dx \right)^{\frac{p(1 - \theta_q)}{q \theta_q}}\; .
\]

\n Expanding the metric $g$ in normal coordinates around $x_q$, one locally gets

\[
(1 - c d_g(x,x_q)^2) dv_g \leq dx \leq (1 + c d_g(x,x_q)^2) dv_g
\]

\n and

\[
|\nabla(u_q \eta_q)|^p \leq |\nabla_g(u_q \eta_q)|^p (1 + c \; d_g(x,x_q)^2) \; .
\]

\n Thus,

\[
\left( \int_{B(0,2)} u_q^p \eta_q^p\; dx \right)^{\frac{1}{\theta_q}} \leq \left( A_0(p,q,p) A_q \int_{B(x_q,2)} |\nabla_g (u_q \eta_q)|^p\; dv_g + c A_q \int_{B(x_q,2)} |\nabla_g (u_q\eta_q)|^p d_g(x,x_q)^2\; dv_g \right)
\]

\[
\times \left( \frac{\int_{B(0,2)} u_q^q \eta_q^q\; dx}{\int_M u_q^q\; dv_g} \right)^{\frac{p(1 - \theta_q)}{q \theta_q}} \; .
\]

\n Using now the inequality

\[
|\nabla_g (u_q \eta_q)|^p \leq |\nabla_g u_q|^p \eta_q^p + c |\eta_q \nabla_g u_q|^{p - 1} |u_q \nabla_g \eta_q| + c |u_q \nabla_g \eta_q|^p \; ,
\]

\n and denoting

\[
X_q =  A_q \int_M \eta_q^p |\nabla_g u_q|^p d_g(x,x_q)^2\; dv_g
\]

\n and

\[
Y_q = A_q \int_M |\nabla_g u_q|^{p -1} |\nabla_g \eta_q| u_q\; dv_g \; ,
\]

\n we deduce that

\begin{equation}\label{af1}
\left( \int_{B(0,2)} u_q^p \eta_q^p\; dx \right)^{\frac{1}{\theta_q}} \leq \left( A_0(p,q,p) A_q \int_M |\nabla_g u_q|^p \eta_q^p\; dv_g + c X_q + c Y_q + c A_q \right)
\left( \frac{\int_{B(0,2)} u_q^q \eta_q^q\; dx}{\int_M u_q^q\; dv_g} \right)^{\frac{p(1 - \theta_q)}{q \theta_q}}\; .
\end{equation}

\n On the other hand, choosing $u_q \eta_q^p$ as a test function in (\ref{l3ep}) and using (\ref{nu}), one gets

\[
A_0(p,q,p) A_q \int_M |\nabla_g u_q|^p \eta_q^p\; dv_g \leq 1 - A_0(p,q,p) A_q C_q + \frac{1}{\theta_q} \left( \int_M u_q^p \eta_q^p\; dv_g - \frac{\int_M u_q^q \eta_q^p\; dv_g}{\int_M u_q^q\;
dv_g} \right)
\]

\[
+ c \int_M |\nabla_g u_q|^{p -1} |\nabla_g \eta_q| u_q\; dv_g \; .
\]

\n Using H\"{o}lder's inequality, (\ref{ldf}), Proposition \ref{P.1} and Lemma \ref{UE} with a suitable value of $\lambda$, the last integral can be estimated as

\[
\int_M |\nabla_g u_q|^{p -1} |\nabla_g \eta_q| u_q\; dv_g \leq c \left( \int_M |\nabla_g u_q|^p\; dv_g \right)^{\frac{p-1}{p}} \left( \int_{B(x_q,2) \setminus B(x_q,1)} u_q^p\; dv_g \right)^{\frac{1}{p}} \leq c A_q\; ,
\]

\n so that

\begin{equation}\label{af2}
A_0(p,q,p) A_q \int_M |\nabla_g u_q|^p \eta_q^p\; dv_g \leq 1 - c_1 A_q C_q + \frac{1}{\theta_q} \left( \int_M u_q^p \eta_q^p\; dv_g - \frac{\int_M u_q^q \eta_q^p\; dv_g}{\int_M u_q^q\;
dv_g} \right) + c A_q
\end{equation}

\n for $q$ close to $p$.

\n Let

\[
Z_q = \frac{1}{\theta_q} \left( \int_M u_q^p \eta_q^p\; dv_g - \frac{\int_M u_q^q \eta_q^q\; dv_g}{\int_M u_q^q\; dv_g} \right) \; .
\]

\n By the mean value theorem and Lemma \ref{UE}, there exists $\gamma_q \in (q,p)$ such that

\begin{equation}\label{af3}
\left| Z_q - \frac{1}{\theta_q}\left( \int_M u_q^p \eta_q^p\; dv_g - \frac{\int_M u_q^q \eta_q^p\; dv_g}{\int_M u_q^q\; dv_g} \right)
\right| \leq \frac{1}{\theta_q} \left| \frac{\int_M u_q^q (\eta_q^q - \eta_q^p)\; dv_g}{\int_M u_q^q\; dv_g} \right| \leq \frac{pq - nq + np}{n} \frac{\int_M | \log \eta_q |
\eta_q^{\gamma_q} u_q^q\; dv_g}{\int_M u_q^q\; dv_g}
\end{equation}

\[
\leq c \frac{\int_{B(x_q,2) \setminus B(x_q,1)} u_q^q\; dv_g}{\int_M u_q^q\; dv_g} \leq c A_q
\]

\n for $q$ close to $p$.

\n Plugging (\ref{af3}) into (\ref{af2}) and after this one into (\ref{af1}), one arrives at

\begin{equation} \label{af4}
\left( \int_{B(0,2)} u_q^p \eta_q^p\; dx \right)^{\frac{1}{\theta_q}} \leq \left( 1 - c A_q C_q + Z_q + c X_q + c Y_q + c A_q\right) \left( \frac{\int_{B(0,2)} u_q^q \eta_q^q\; dx}{\int_M u_q^q\; dv_g} \right)^{\frac{p(1 - \theta_q)}{q \theta_q}}
\end{equation}

\n for $q$ close to $p$.

\n In order to estimate $X_q$ and $Y_q$, we take $u_q d_g^2 \eta_q^p$ as a test function in (\ref{l3ep}). From this choice, we derive

\[
X_q \leq \frac{\nu_q}{\theta_q}\left( \int_M u_q^p \eta_q^p d_g(x,x_q)^2\; dv_g - \frac{\int_M u_q^q \eta_q^q d_g(x,x_q)^2\; dv_g}{\int_M u_q^q\;
dv_g}\right) + c A_q \int_M u_q \eta_q^p |\nabla_g u_q|^{p - 1} d_g(x,x_q)\; dv_g
\]

\[
+ c \frac{\int_M u_q^q \eta_q^q d_g(x,x_q)^2\; dv_g}{\int_M u_q^q dv_g} + c Y_q + c A_q\; .
\]

\n We now estimate the first two terms of the right-hand side above. Namely, after a change of variable, one has

\[
\frac{1}{\theta_q} \int_{M} \left| u_q^p \eta_q^p d_g(x,x_q)^2 - \frac{u_q^q \eta_q^q d_g(x,x_q)^2}{\int_M u_q^q\; dv_g}\right|\; dv_g \leq c A_q^{\frac{2}{p}} \frac{1}{\theta_q} \int_{B(0,2 A_q^{- \frac{1}{p}})} \left| \varphi_q^p \tilde{\eta}_q^p - \varphi_q^q \tilde{\eta}_q^q \right| |x|^2\; dx
\]

\[
= c A_q^{\frac{2}{p}} \int_{B(0,2 A_q^{- \frac{1}{p}})} \left| \log(\varphi_q \tilde{\eta}_q) \right| (\varphi_q \tilde{\eta}_q)^{\rho_q} |x|^2\; dx
\]

\n for some $\rho_q \in (q,p)$, where $\tilde{\eta}_q(x) = \eta_q(A_q^{\frac{1}{p}} x)$.

\n Thus, using Lemma \ref{UE} and the assumption $p \leq 2$, one obtains

\begin{equation} \label{est6}
\frac{1}{\theta_q} \int_{M} \left| u_q^p \eta_q^p d_g(x,x_q)^2 - \frac{u_q^q \eta_q^q d_g(x,x_q)^2}{\int_M u_q^q\; dv_g}\right|\; dv_g \leq c A_q\; .
\end{equation}

\n In particular, since

\[
\int_M u_q^p \eta_q^p d_g(x,x_q)^2\; dv_g - \frac{\int_M u_q^q \eta_q^q d_g(x,x_q)^2\; dv_g}{\int_M u_q^q\; dv_g} = \int_{M} \left( u_q^p \eta_q^p d_g(x,x_q)^2 - \frac{u_q^q \eta_q^q d_g(x,x_q)^2}{\int_M u_q^q \; dv_g}\; \right) dv_g\; ,
\]

\n we have

\[
\frac{1}{\theta_q}\left| \int_M u_q^p \eta_q^p d_g(x,x_q)^2\; dv_g - \frac{\int_M u_q^q \eta_q^q d_g(x,x_q)^2\; dv_g}{\int_M u_q^q\;
dv_g}\right| \leq c A_q
\]

\n for $q$ close to $p$.

\n Using H\"{o}lder's inequality, (\ref{ldf}), Proposition \ref{P.1}, Lemma \ref{UE} and the fact that $p \leq 2$, one gets

\[
\int_M u_q \eta_q^p |\nabla_g u_q|^{p - 1} d_g(x,x_q)\; dv_g \leq c \left( \int_M |\nabla_g u_q|^p\; dv_g \right)^{\frac{p-1}{p}} \left( \int_{B(x_q,2)} u_q^p d_g(x,x_q)^p\; dv_g \right)^{\frac{1}{p}} \]

\[
\leq c A_q^{\frac{2-p}{p}} \left( \int_{B(0,2 A_q^{- \frac{1}{p}})} \varphi_q^p |x|^p\; dx \right)^{\frac{p-1}{p}} \leq c
\]

\n for $q$ close to $p$.

\n Consequently, the above estimates guarantees that

\begin{equation}\label{f2}
X_q \leq  c \frac{\int_M u_q^q \eta_q^q d_g(x,x_q)^2\; dv_g}{\int_M u_q^q\; dv_g} + c Y_q + c A_q\; .
\end{equation}

\n Again, with the aid of Lemma \ref{UE} and the condition $p \leq 2$, one derives

\begin{equation}\label{f1}
\frac{\int_M u_q^q \eta_q^q d_g(x,x_q)^2\; dv_g}{\int_M u_q^q\; dv_g} \leq c A_q^{\frac{2}{p}}  \int_{B(0,2A_q^{-1/p})} \varphi_q^q |x|^2\; dh_q \leq c A_q
\end{equation}

\n By H\"{o}lder's inequality, (\ref{ldf}) and Proposition \ref{P.1}, we also have

\begin{equation}\label{est3}
Y_q \leq c A_q^{\frac{1}{p}} \int_{B(x_q,2) \backslash B(x_q,1)} u_q^p\; dv_g \leq c A_q\; .
\end{equation}

\n In conclusion, from (\ref{f2}), (\ref{f1}) and (\ref{est3}), it follows that

\[
\left(\int_{B(0,2)} u_q^p \eta_q^p\; dx \right)^{\frac{1}{\theta_q}} \leq \left(1 + Z_q - c_1 A_q C_q + c A_q \right) \left(\frac{\int_{B(0,2)} (u_q \eta_q)^q\; dx}{\int_M u_q^q\; dv_g} \right)^{\frac{p(1 - \theta_q)}{q \theta_q}}
\]

\n for $q$ close to $p$.

\n Taking logarithm of both sides and using the fact that $\frac{p(1 - \theta_q)}{q \theta_q} = \frac{1}{\theta_q} - \frac{n - p}{n}$, one has

\begin{equation}\label{d1}
\frac{1}{\theta_q} \left( \log \int_{B(0,2)} u_q^p \eta_q^p dx - \log \left( \frac{\int_{B(0,2)} u_q^q \eta_q^q dx}{\int_M u_q^q\; dv_g}
\right) \right) \leq \log (1 + Z_q - c_1 A_q C_q + c A_q) - \frac{n - p}{n} \log \left( \frac{\int_{B(0,2)} u_q^q \eta_q^q\; dx}{\int_M u_q^q\; dv_g}
\right)\; .
\end{equation}

\n By the mean value theorem,

\begin{equation}\label{d2}
\log \int_{B(0,2)} u_q^p \eta_q^p\; dx - \log \left( \frac{\int_{B(0,2)} u_q^q \eta_q^q\; dx}{\int_M u_q^q\; dv_g} \right)
= \frac{1}{\tau_q} \left( \int_{B(0,2)} u_q^p \eta_q^p\; dx - \frac{\int_{B(0,2)} u_q^q \eta_q^q\; dx}{\int_M u_q^q\; dv_g} \right)
\end{equation}

\n for some number $\tau_q$ between the expressions

\[
\int_{B(0,2)} u_q^p \eta_q^p\; dx\ \ {\rm and}\ \ \frac{\int_{B(0,2)} u_q^q \eta_q^q\; dx}{\int_M u_q^q\; dv_g} \; .
\]

\n Using Cartan's expansion of $g$ in normal coordinates around $x_q$ and Lemma \ref{UE}, one obtains

\begin{equation} \label{est4}
\max\left\{ \left| \int_{B(0,2)} u_q^p \eta_q^p\; dx - \int_M u_q^p \eta_q^p\; dv_g \right|, \left| \frac{\int_{B(0,2)} u_q^q \eta_q^q\; dx}{\int_M u_q^q\; dv_g} - \frac{\int_M u_q^q \eta_q^q dv_g}{\int_M u_q^q dv_g} \right|\right\} \leq c A_q
\end{equation}

\n for $q$ close to $p$.

\n Indeed, since $p \leq 2$,

\[
\left| \int_{B(0,2)} u_q^p \eta_q^p\; dx - \int_M u_q^p \eta_q^p\; dv_g \right| \leq c \int_M u_q^p \eta_q^p d_g(x,x_q)^2 \; dv_g \leq c A_q^{\frac{2}{p}}  \int_{B(0,2A_q^{-1/p})} \varphi_q^p |x|^2\; dh_q \leq c A_q
\]

\n and, by (\ref{f1}),

\[
\left| \frac{\int_{B(0,2)} u_q^q \eta_q^q\; dx}{\int_M u_q^q\; dv_g} - \frac{\int_M u_q^q \eta_q^q dv_g}{\int_M u_q^q dv_g} \right| \leq c \frac{\int_M u_q^q \eta_q^q d_g(x,x_q)^2\; dv_g}{\int_M u_q^q\; dv_g} \leq c A_q\; .
\]

\n Moreover, we also have

\begin{equation} \label{est5}
\max\left\{ \left| \int_M (u_q \eta_q)^p\; dv_g - 1\right|, \left|\frac{\int_M u_q^q \eta_q^q\; dv_g}{\int_M u_q^q\; dv_g} - 1 \right|\right\} \leq c A_q
\end{equation}

\n for $q$ close to $p$.

\n In fact, by Lemma \ref{UE},

\[
\left| \int_M u_q^p \eta_q^p\; dv_g - 1\right| = \left| \int_M u_q^p \eta_q^p\; dv_g - \int_M u_q^p \; dv_g\right| \leq c \int_{M \setminus B(x_q,1)} u_q^p \; dv_g \leq c A_q
\]

\n and

\[
\left|\frac{\int_M u_q^q \eta_q^q\; dv_g}{\int_M u_q^q\; dv_g} - 1 \right| \leq c \frac{\int_{M \setminus B(x_q,1)} u_q^q\; dv_g}{\int_M u_q^q\; dv_g} \leq c A_q\; .
\]

\n Thanks to (\ref{est4}) and (\ref{est5}), one easily deduces that $\tau_q^{-1} = 1 + O(A_q)$. So, by (\ref{d2}),

\begin{equation}\label{d3}
\frac{1}{\theta_q} \left( \log \int_{B(0,2)} u_q^p \eta_q^p\; dx - \log \left( \frac{\int_{B(0,2)} u_q^q \eta_q^q\; dx}{\int_M u_q^q\; dv_g} \right) \right)
= \frac{1}{\theta_q} \left( \int_{B(0,2)} u_q^p \eta_q^p\; dx - \frac{\int_{B(0,2)} u_q^q \eta_q^q\; dx}{\int_M u_q^q\; dv_g} \right) (1 + O(A_q))\; .
\end{equation}

\n Now, from Cartan's expansion and (\ref{est6}), note that

\[
\frac{1}{\theta_q} \left( \int_{B(0,2)} u_q^p \eta_q^p\; dx - \frac{\int_{B(0,2)} u_q^q \eta_q^q\; dx}{\int_M u_q^q\; dv_g} \right) = \frac{1}{\theta_q} \left(\int_{B(0,2)} u_q^p \eta_q^p - \frac{u_q^q \eta_q^q}{\int_M u_q^q\; dv_g}\; dx \right)
\]

\[
=  \frac{1}{\theta_q} \left(\int_M u_q^p \eta_q^p - \frac{u_q^q \eta_q^q}{\int_M u_q^q\; dv_g}\; dv_g \right) + \frac{1}{\theta_q} \left(\int_M \left( u_q^p \eta_q^p - \frac{u_q^q \eta_q^q}{\int_M u_q^q\; dv_g} \right) O(d_g(x, x_q)^2)\; dv_g \right)
\]

\[
= Z_q + O(A_q)
\]

\n for $q$ close to $p$.

\n Replacing this inequality in (\ref{d3}), one gets

\[
\frac{1}{\theta_q} \left(\log \int_{B(0,2)} u_q^p \eta_q^p\; dx - \log \left( \frac{\int_{B(0,2)} u_q^q \eta_q^q\; dx}{\int_M u_q^q\; dv_g} \right)\right) \geq Z_q - c A_q\; .
\]

\n In turn, plugging the above inequality in (\ref{d1}), one has

\[
Z_q - c A_q \leq \log (1 + Z_q - c_1 A_q C_q + c A_q) + \frac{n - p}{n} \left| \log \left( \frac{\int_{B(0,2)} u_q^q \eta_q^q\; dx}{\int_M u_q^q\; dv_g} \right) \right|\; .
\]

\n Finally, using Cartan's expansion of $g$ in normal coordinates, Taylor's expansion of the function $\log$ and Lemma \ref{UE}, one obtains

\[
\left| \log \left( \frac{\int_{B(0,2)} u_q^q \eta_q^q\; dx}{\int_M u_q^q\; dv_g} \right) \right| \leq c \frac{\int_{M \backslash B(x_q,1)} u_q^q\; dv_g}{\int_M u_q^q\; dv_g} + c \frac{\int_M u_q^q \eta_q^q d_g(x,x_q)^2\; dv_g}{\int_M u_q^q\; dv_g} \leq c A_q\; .
\]

\n In short, for each fixed $c > 0$ large, we deduce that

\[
Z_q \leq \log (1 + Z_q - c_1 A_q C_q + c A_q) + c A_q \; .
\]

\n Since $C_q \rightarrow + \infty$, there exists $q_0 = q_0(c)$ such that

\[
1 + Z_q - c_1 A_q C_q + c A_q \leq Z_q + 1 - c A_q \leq \log (1 + Z_q - c_1 A_q C_q + c A_q) + 1
\]

\n for all $q_0 < q < p$. Note that the preceding inequality implies $1 + Z_q - c_1 A_q C_q + c A_q = 1$ for all $q > q_0$. But this is a contradiction since $c$ can be chosen large enough. \bl \\

{\bf Acknowledgements.} The first author was partially supported by CAPES through INCTmat and the second one was partially supported by CNPq and Fapemig.

\end{document}